\documentclass[12pt,english,british,refpage,intoc,bibliography=totoc,index=totoc,BCOR=7.5mm,captions=tableheading]{extarticle}
\usepackage[latin9]{inputenc}
\usepackage[a4paper]{geometry}
\usepackage{color}
\usepackage{babel}
\usepackage{mathrsfs}
\usepackage{amsmath}
\usepackage{amsthm}
\usepackage{amssymb}
\usepackage[numbers]{natbib}
\usepackage[unicode=true,pdfusetitle,
 bookmarks=true,bookmarksnumbered=true,bookmarksopen=false,
 breaklinks=false,pdfborder={0 0 1},backref=false,colorlinks=true]
 {hyperref}
\hypersetup{
 pdfborderstyle=,linkcolor=black,citecolor=blue,urlcolor=black,filecolor=blue,pdfpagelayout=OneColumn,pdfnewwindow=true,pdfstartview=XYZ,plainpages=false}

\makeatletter
\theoremstyle{plain}
\newtheorem{thm}{\protect\theoremname}[section]
\theoremstyle{plain}
\newtheorem{lem}[thm]{\protect\lemmaname}
\theoremstyle{definition}
\newtheorem{example}[thm]{\protect\examplename}
\theoremstyle{plain}
\newtheorem{cor}[thm]{\protect\corollaryname}

\@ifundefined{date}{}{\date{}}
\usepackage{babel}
\usepackage{caption}
\usepackage[nottoc]{tocbibind}
\allowdisplaybreaks[4]
\usepackage{dsfont}
\DeclareMathAlphabet{\mathcal}{OMS}{cmsy}{m}{n}

\addto\captionsbritish{}
\addto\captionsbritish{\renewcommand{\lemmaname}{Lemma}}
\addto\captionsbritish{}
\addto\captionsbritish{\renewcommand{\theoremname}{Theorem}}
\addto\captionsenglish{}
\addto\captionsenglish{\renewcommand{\lemmaname}{Lemma}}
\addto\captionsenglish{}
\addto\captionsenglish{\renewcommand{\theoremname}{Theorem}}

\providecommand{\lemmaname}{Lemma}

\providecommand{\theoremname}{Theorem}

\addto\captionsbritish{}
\addto\captionsbritish{\renewcommand{\lemmaname}{Lemma}}
\addto\captionsbritish{}
\addto\captionsbritish{\renewcommand{\theoremname}{Theorem}}
\addto\captionsenglish{}
\addto\captionsenglish{\renewcommand{\lemmaname}{Lemma}}
\addto\captionsenglish{}
\addto\captionsenglish{\renewcommand{\theoremname}{Theorem}}

\addto\captionsbritish{\renewcommand{\corollaryname}{Corollary}}
\addto\captionsbritish{}
\addto\captionsbritish{\renewcommand{\lemmaname}{Lemma}}
\addto\captionsbritish{\renewcommand{\theoremname}{Theorem}}
\addto\captionsenglish{\renewcommand{\corollaryname}{Corollary}}
\addto\captionsenglish{}
\addto\captionsenglish{\renewcommand{\lemmaname}{Lemma}}
\addto\captionsenglish{\renewcommand{\theoremname}{Theorem}}
\providecommand{\corollaryname}{Corollary}

\makeatother

\addto\captionsbritish{\renewcommand{\corollaryname}{Corollary}}
\addto\captionsbritish{\renewcommand{\examplename}{Example}}
\addto\captionsbritish{\renewcommand{\lemmaname}{Lemma}}
\addto\captionsbritish{\renewcommand{\theoremname}{Theorem}}
\addto\captionsenglish{\renewcommand{\corollaryname}{Corollary}}
\addto\captionsenglish{\renewcommand{\examplename}{Example}}
\addto\captionsenglish{\renewcommand{\lemmaname}{Lemma}}
\addto\captionsenglish{\renewcommand{\theoremname}{Theorem}}
\providecommand{\corollaryname}{Corollary}
\providecommand{\examplename}{Example}
\providecommand{\lemmaname}{Lemma}
\providecommand{\theoremname}{Theorem}

\begin{document}
\global\long\def\Sgm{\boldsymbol{\Sigma}}%
\global\long\def\W{\boldsymbol{W}}%
\global\long\def\H{\boldsymbol{H}}%
\global\long\def\P{\mathbb{P}}%
\global\long\def\Q{\mathbb{Q}}%
\global\long\def\E{\mathbb{E}}%
\global\long\def\R{\mathbb{R}}%
\global\long\def\PP{\mathbb{P}}%
\global\long\def\n{{(n)}}%
\global\long\def\nn{{(n+1)}}%
\global\long\def\nk{{(n_{k})}}%
\global\long\def\xx{\boldsymbol{x}}%
\global\long\def\yy{\boldsymbol{y}}%
\global\long\def\zz{\boldsymbol{z}}%
\global\long\def\xixi{\boldsymbol{\xi}}%
\global\long\def\omegaomega{\boldsymbol{\omega}}%
\global\long\def\ff{\boldsymbol{f}}%
\global\long\def\uu{\boldsymbol{u}}%
\global\long\def\KK{\boldsymbol{K}}%
\global\long\def\E{\mathbb{E}}%
\global\long\def\R{\mathbb{R}}%
\global\long\def\PP{\mathbb{P}}%
\global\long\def\sgn{\mathrm{sgn}}%
\global\long\def\n{{(n)}}%
\global\long\def\nn{{(n+1)}}%
\global\long\def\nj{{(n_{j})}}%
\global\long\def\np{{(n_{p})}}%
\global\long\def\nk{{(n_{k})}}%
\global\long\def\bvert{\big\vert}%
\global\long\def\Bvert{\Big\vert}%
\global\long\def\bbvert{\bigg\vert}%
\global\long\def\bVert{\big\Vert}%
\global\long\def\BVert{\Big\Vert}%
\global\long\def\bbVert{\bigg\Vert}%
\global\long\def\FF{\mathcal{F}}%
\global\long\def\Fscr{\mathscr{F}}%
\global\long\def\eat{e^{\alpha t}}%
\global\long\def\eas{e^{\alpha s}}%
\global\long\def\Ytxi{Y_{t}^{\xi}}%
\global\long\def\Yteta{Y_{t}^{\eta}}%
\global\long\def\Ysxi{Y_{s}^{\xi}}%
\global\long\def\Yseta{Y_{s}^{\eta}}%
\global\long\def\supt{\sup_{0\leq t\leq T}}%
\global\long\def\equalDistrib{\,\frac{\buildrel\Delta}{=}\,}%
\numberwithin{equation}{section} 
\global\long\def\blue#1{\textcolor{blue}{#1}}%
\global\long\def\red#1{\textcolor{red}{#1}}%
\global\long\def\ww{\boldsymbol{w}}%
\global\long\def\xx{\boldsymbol{x}}%
\global\long\def\yy{\boldsymbol{y}}%
\global\long\def\zz{\boldsymbol{z}}%
\global\long\def\bb{\boldsymbol{b}}%
\global\long\def\ff{\boldsymbol{f}}%
\global\long\def\kk{\boldsymbol{k}}%
\global\long\def\rr{\boldsymbol{r}}%
\global\long\def\xixi{\boldsymbol{\xi}}%
\global\long\def\omegaomega{\boldsymbol{\omega}}%
\global\long\def\alphaalpha{\boldsymbol{\alpha}}%
\global\long\def\betabeta{\boldsymbol{\beta}}%
\global\long\def\gammagamma{\boldsymbol{\gamma}}%
\global\long\def\kappakappa{\boldsymbol{\kappa}}%
\global\long\def\mumu{\boldsymbol{\mu}}%
\global\long\def\thetatheta{\boldsymbol{\theta}}%
\global\long\def\lambdalambda{\boldsymbol{\lambda}}%
\global\long\def\uu{\boldsymbol{u}}%
\global\long\def\ii{\boldsymbol{i}}%
\global\long\def\jj{\boldsymbol{j}}%
\global\long\def\boldg{\boldsymbol{g}}%
\global\long\def\KK{\boldsymbol{K}}%
\global\long\def\GG{\boldsymbol{G}}%
\global\long\def\MM{\mathscr{M}}%
\global\long\def\Osc{\text{Osc}}%
\global\long\def\tildeX{\widetilde{X}}%
\global\long\def\tildeu{\widetilde{u}}%

\title{McKean-Vlasov type stochastic differential equations arising from
the random vortex method}
\author{Zhongmin Qian\thanks{Mathematical Institute, University of Oxford, Oxford OX2 6GG, and
Oxford Suzhou Centre for Advanced Research. \emph{Email}: $\mathtt{qianz@maths.ox.ac.uk}$} $\;$and Yuhan Yao\thanks{Mathematical Institute, University of Oxford, Oxford OX2 6GG. \emph{Email}:
$\mathtt{yaoy@maths.ox.ac.uk}$}}
\maketitle
\begin{abstract}
We study a class of McKean-Vlasov type stochastic differential equations
(SDEs) which arise from the random vortex dynamics and other physics
models. By introducing a new approach we resolve the existence and
uniqueness of both the weak and strong solutions for the McKean-Vlasov
stochastic differential equations whose coefficients are defined in
terms of singular integral kernels such as the Biot-Savart kernel.
These SDEs which involve the distributions of solutions are in general
not Lipschitz continuous with respect to the usual distances on the
space of distributions such as the Wasserstein distance. Therefore
there is an obstacle in adapting the ordinary SDE method for the study
of this class of SDEs, and the conventional methods seem not appropriate
for dealing with such distributional SDEs which appear in applications
such as fluid mechanics.

\selectlanguage{english}%
\medskip{}

\selectlanguage{british}%
\emph{key words}: Aronson estimates, Cameron-Martin formula, diffusion
processes, McKean-Vlasov SDEs, strong solution, vorticity
equation

\medskip{}

\emph{MSC classifications}: 60H30, 35Q30, 35Q35, 76D03, 76D05, 76D17
\end{abstract}

\section{Introduction}

In this paper, we study the following McKean-Vlasov type stochastic
differential equations
\begin{equation}
\begin{cases}
dX^{i}(x,t)=\left(\sum_{j=1}^{d}\int_{\mathbb{R}^{d}}\left.\mathbb{E}\left[K_{j}^{i}(z-X(y,t))\right]\right|_{z=X(x,t)}\omega_{0}^{j}(y)dy\right)dt+\sqrt{2\nu}dB^{i}(t),\\
X(x,t)=x,\textrm{ for }x\in\mathbb{R}^{d}
\end{cases}\label{eq:Mckean-Vlasov}
\end{equation}
where $i=1,\ldots,d$, $\nu>0$ is a constant (which has its origin
in fluid mechanics, namely the kinetic viscosity), $B=(B^{1},\ldots,B^{d})$
is a standard Brownian motion on some probability space and $\omega_{0}=(\omega_{0}^{1},\ldots,\omega_{0}^{d})$
is the initial data to the corresponding non-linear (and non-local)
partial differential equations (PDEs), see \eqref{eq:PDE1} below.
The structure kernel function which defines SDE \eqref{eq:Mckean-Vlasov}
$K=(K_{j}^{i})$ is a $d\times d$ matrix-valued Borel measurable
function on $\mathbb{R}^{d}$, which is continuous except at several
singularities. The study of  SDE \eqref{eq:Mckean-Vlasov} is inspired from the
random vortex method in fluid mechanics, in which the integral kernel
$K$ is singular at $0$. We will explain the random vortex model
and formulate the problem \eqref{eq:Mckean-Vlasov} more precisely
in the next section. 

\eqref{eq:Mckean-Vlasov} is a system of SDEs which involves the distributions
of its solutions. This type of SDEs and SDEs which share the same
nature may arise from physics models and from applied mathematics, and
they have been studied intensively over the past decades. There is
a large amount of literature devoting to various aspects of McKean-Vlasov
equations, initiated by McKean in his seminal paper \citep{McKea 1966}
(see for example \citep{Cameron Delarue Vol1,Funaki 1984,Sznitman 1991}
for some recent progress, \citep{buckdahn2009mean,buckdahn2017mean,bauer2018strong,Meleard 2000,raynal2020strong,Tomasevic 2021,Talay Tomasevic 2020}
and the literature therein).

The study of McKean-Valsov SDEs and the renewed interest in SDEs involving
solution distributions in recent years are largely influenced by their
connections with some non-local and non-linear PDEs arising from physics
models. In this aspect, McKean-Valsov type SDEs provide the theoretical
foundation for numerical methods such as the particle method for simulating
the solutions to this kind of PDEs. For example, the propagation of
chaos (law of large numbers) for solving the corresponding PDE of
\eqref{eq:Mckean-Vlasov} may be formulated by replacing the expectation
by the empirical measure, to obtain the following system
\begin{equation}
\begin{cases}
dX^{n,\boldsymbol{k}}=\frac{1}{N}\sum_{n=1}^{N}\sum_{\boldsymbol{j}}\varepsilon^{d}K(X^{n,\boldsymbol{k}}-X^{n,\boldsymbol{j}})\omega^{\boldsymbol{j}}dt+\sqrt{2\nu}dB^{n}(t),\\
X^{n,\boldsymbol{k}}=\varepsilon\boldsymbol{k},\textrm{ for }\boldsymbol{k}\in\mathbb{Z}^{d},
\end{cases}\label{eq:sN}
\end{equation}
where $B^{n}$ are independent copies of $d$-dimensional Brownian
motion, and $\varepsilon>0$ is the lattice size. The previous random
system is the essential ingredient in the random vortex method, see
for example \citep{marchioro1982hyrodynamics,Meleard 2000,osada1987propagation}.
Other numerical approximations have also been employed to look for
large deviation  results, see for example \citep{beale1982vortex,cottet1988new,long1988convergence,majda2002vorticity}
for detailed discussions. The difficulty however, in particular in
the case that the dimension $d=3$ and $K$ is the Biot-Savart kernel,
comes from the fact that the kernel $K$ is too singular at $0$,
hence the Lipschitz continuity of the coefficients appearing in \eqref{eq:Mckean-Vlasov},
which is essential  (see for example \citep{Cameron Delarue Vol1,Sznitman 1991}),
can not be expected. 

SDE \eqref{eq:Mckean-Vlasov} has an independent interest by its own
of course besides its significance in fluid dynamics. The research
for this type of SDEs has been dominated, to the best knowledge of
the present authors, by the use of It\^o's SDE theory in one or another
way which requires the Lipschitz continuity of $K$ with respect to
the variational distance or the Wasserstein distance when one seeks
for strong solutions, or by means of martingale problem for weak solutions.
Unfortunately, these approaches are not appropriate for the study
of \eqref{eq:Mckean-Vlasov} when $K$ is singular such as the Biot-Savart
kernel $-\frac{1}{4\pi}\frac{x}{|x|^{3}}$ (where $d=3$) which is
explored near zero like $1/|x|^{2}$.

In the present paper, we overcome these difficulties by devising a
new and powerful approach which allows us to establish the existence
and uniqueness of strong and weak solutions of \eqref{eq:Mckean-Vlasov}
under very weak conditions on the singular integral kernel $K$. In
particular, our results apply to the the Biot-Savart kernel in any
dimension, and also apply to the Green kernels (such as $\ln|x|$
in dimension $2$, $1/|x|^{d-2}$ for $d>2$), the Riesz kernels $1/|x|^{\gamma}$
where $\gamma\in[0,d)$ on $\mathbb{R}^{d}$ and many other singular
integral kernels.

Our novel approach is based on the following simple observation. If
$K$ is singular, then the mapping
\[
(x,\mu)\rightarrow\sum_{j=1}^{d}\int_{\mathbb{R}^{d}}\mathbb{E}\left[K_{j}^{i}(x-\xi)\right]\omega_{0}^{j}(y)dy
\]
where $\xi$ has a distribution $\mu$, is unlikely Lipschitz continuous
with respect to the variational or the Wasserstein metric on the space
of distributions. However we recognise that the distributions of possible
solutions to \eqref{eq:Mckean-Vlasov}, even $K$ is singular, have
much higher regularity. In fact, if $\{X(x,t):x\in\mathbb{R}^{d},t\geq0\}$
is a (weak) solution to \eqref{eq:Mckean-Vlasov} then
\begin{equation}
b^{i}(x,t)=\sum_{j=1}^{d}\int_{\mathbb{R}^{d}}\mathbb{E}\left[K_{j}^{i}(x-X(y,t))\right]\omega_{0}^{j}(y)dy\label{eq:drift}
\end{equation}
defines a vector field (although the vector field $b(x,t)$ is defined
via the solution of the SDE), and $X(x,t)$ must be a weak solution
to the diffusion process defined by ordinary SDE
\[
dX_{t}=b(X_{t},t)dt+\sqrt{2\nu}dB_{t}.
\]
Therefore the distribution of $X(x,t)$ can be represented by Cameron-Martin
formula in terms of the Wiener measure and many results from diffusion
processes thus can be brought in to the study of SDE \eqref{eq:Mckean-Vlasov}.
In this paper the major technical tool is the sharp heat kernel estimates
obtained in \citep{qian2002sharp,qian2003comparison}. 

The paper is organised as the following. In Section 2, we recall the
random vortex problem and derive the SDE \eqref{eq:Mckean-Vlasov}
from the vorticity equation and formulate SDE \eqref{eq:Mckean-Vlasov}
in a form which will be appropriate in the work frame of the present
paper. In Section 3, we collect a few facts about diffusion processes
whose infinitesimal generators are elliptic operators of second order,
and we prove several technical estimates which will be used to prove
our main results. In Section 4, we define a non-linear mapping associated
with SDE \eqref{eq:Mckean-Vlasov} and prove it is a contractive mapping,
then we show that \eqref{eq:Mckean-Vlasov} admits a unique weak solution.
In Section 5 we show that a strong solution can be constructed, and
show that the drift vector field \eqref{eq:drift} is H\"older continuous.
Section 6 recovers solutions to the corresponding non-linear PDEs
by using the solutions to \eqref{eq:Mckean-Vlasov}, which can be
considered as a probabilistic representation for this class of non-local
and non-linear PDEs.

\emph{Convention on Notations}. The following set of conventions are
employed throughout the paper. Firstly Einstein's convention on summation
on repeated indices through their ranges is assumed, unless otherwise
specified. If $A$ is a vector or a vector field (in the space of
dimension $d$) dependent on some parameters, then its components
are labelled with upper-script indices so that $A=\left(A^{i}\right)=\left(A^{1},\ldots,A^{d}\right)$.
The same convention applies to coordinates too, so that $x=(x^{i})=(x^{1},\ldots,x^{d})$.
If $u$ is a vector field on $\mathbb{R}^{3}$ then $\nabla\wedge u$
denotes the curl of $u$ which is again a vector field on $\mathbb{R}^{3}$
with its component $\varepsilon^{ijk}\frac{\partial u^{k}}{\partial x^{j}}$.
If $f$ is a function on $\mathbb{R}^{d}$, then $\left\Vert f\right\Vert _{p}$
or $\left\Vert f\right\Vert _{L^{p}(\mathbb{R}^{d})}$ denotes the
$L^{p}$-norm of $f$ with respect to the Lebesgue measure. Similarly,
if $f(x,t)$ is a function on $\mathbb{R}^{d}\times[0,T]$ then $\left\Vert f\right\Vert _{L^{p}(\mathbb{R}^{d}\times[0,T])}$
or, if no confusion is possible, $\left\Vert f\right\Vert _{p}$ denotes
the $L^{p}$-norm with respect to the Lebesgue measure on the product
space $\mathbb{R}^{d}\times[0,T]$.

\section{Random vortex method -- from PDE to SDE}

Particle formulations for fluid flows have been studied as a tool
for understanding fluid dynamics of turbulence. The underlying idea
is simple, originally due to Taylor \citep{Taylor 1921}. Instead
of considering the velocity vector field $u(x,t)$ of the flow, one
may study the dynamics of trajectories $X(x,t)$ of the fluid particles
emitting from $x$ at the moment $0$, i.e. the dynamical equation
\begin{equation}
\frac{d}{dt}X(x,t)=u(X(x,t),t),\quad X(x,0)=x\label{eq:X-1}
\end{equation}
and reformulate the equation of motion of the vorticity $\omega=\nabla\wedge u$
into an evolution equation for $X(x,t)$. This approach works well
for certain inviscid fluids.

For viscous incompressible fluid with constant viscosity $\nu>0$,
a natural idea is to consider Brownian particles instead, i.e. $X(x,t)$
is modelled by the Taylor diffusion 
\[
dX(x,t)=u(X(x,t),t)+\sqrt{2\nu}dB_{t},\quad X(x,0)=x
\]
where $B$ is a standard Brownian motion, and rewrite the equation
of vorticity motion in terms of the distribution of the Taylor diffusion.
This approach is called the random vortex method, see for example
\citep{chorin1973numerical,chorin1994vorticity,cottet2000vortex,majda2002vorticity,saffman1992vortex,ting2007vortex}
and etc. for a comprehensive account including the recent exciting
progress. For incompressible fluid flows, $u(x,t)$ satisfies the
Navier-Stokes equations 
\begin{align}
\frac{\partial}{\partial t}u+u\cdot\nabla u & =\nu\Delta u-\nabla p,\label{eq:nav-1}\\
\nabla\cdot u & =0\label{eq:nav2}
\end{align}
where $p(x,t)$ is a scalar function representing the pressure at
$(x,t)$. If the fluid is constrained in a finite region, then certain
boundary conditions must be identified, but for simplicity, we consider
the case where the evolution of the fluid can take place without physical
boundary and also without external force. In this case, the implicit
boundary condition at infinity is applied: both $u(x,t)$ and $p(x,t)$
tend to zero sufficiently fast as $|x|\rightarrow\infty$. This is
the model used in the homogeneous turbulence for example. The incompressible
condition \eqref{eq:nav2} allows to reformulate the first equation
\eqref{eq:nav-1} in terms of the fluid vorticity $\omega=\nabla\wedge u$
and the equation of vorticity motion is the following vorticity equation
\begin{equation}
\frac{\partial}{\partial t}\omega+u\cdot\nabla\omega=\nu\Delta\omega+\omega\cdot\nabla u,\label{eq:vort-ns1}
\end{equation}
where the velocity field $u$ can be recovered from the Laplace equation
\begin{equation}
\Delta u=-\nabla\wedge\omega.\label{eq:vort-ns2}
\end{equation}

\subsection{Taylor's diffusions}

In our approach, Taylor's diffusions will play a crucial r\^ole,
so the goal of this part is not only for the propose of describing
the vortex dynamics, but also for establishing a few notions and notations
which will be used frequently throughout the paper. The vorticity
equation may be written as 
\begin{equation}
\left(\frac{\partial}{\partial t}-L_{-u}\right)\omega=\omega\cdot\nabla u,\label{eq:dyn-variablew}
\end{equation}
where we have introduced the following notation: if $b(x,t)$ is a
time-dependent vector field (here $t$ is the time variable), then
\begin{equation}
L_{b(x,t)}=\nu\Delta+b(x,t)\cdot\nabla\label{eq:s-new4}
\end{equation}
which is a differential operator of second order and is time inhomogeneous
in general. This convention will be applied to any time dependent
vector field $b(x,t)$ on $\mathbb{R}^{d}$ where $d$ is not necessary
to be $3$. If no confusion may arise, the argument $(x,t)$ will
be suppressed. $L_{b}$ is the infinitesimal generator of the Taylor
diffusion describing the motion of Brownian particles $(X_{t})_{t\geq0}$,
which can be defined by the It\^o's stochastic differential equation
\begin{equation}
dX_{t}=b(X_{t},t)dt+\sqrt{2\nu}dB_{t}\label{eq:s new 5}
\end{equation}
where $(B)_{t\geq0}$ is a standard Brownian motion in $\mathbb{R}^{d}$
on some probability space $(\varOmega,\mathcal{F},\mathbb{P})$. The
transition probability density function of the $L_{b}$-diffusion
is denoted by $p_{b}(\tau,x,t,y)$ for $t>\tau\geq0$ and $x,y\in\mathbb{R}^{d}$
in the sense that
\[
\mathbb{P}\left[X_{t}\in dy|X_{\tau}=x\right]=p_{b}(\tau,x,t,y)dy.
\]

The formal adjoint operator of $L_{b}$ is given by 
\begin{equation}
L_{b}^{\star}=\nu\Delta-b\cdot\nabla-\nabla\cdot b=L_{-b}-\nabla\cdot b\label{eq:L_b star}
\end{equation}
which is again a diffusion operator if and only if the vector field
$b$ is divergence-free. In particular, if $b(\cdot,t)$ is solenoidal
then $L_{b}^{\star}=L_{-b}$. The following lemma contains the facts
about the elliptic operator $L_{b}$ which will be used throughout
the paper.
\begin{lem}
\label{lem:fund-sol-vortex} Let $b(x,t)$ be a time-dependent Borel
measurable and bounded vector field on $\mathbb{R}^{d}$. Let $\varGamma_{F,L_{b}^{\star}}$
be the fundamental solution of the forward parabolic equation $\left(\frac{\partial}{\partial t}-L_{b}^{\star}\right)u=0$
(see \citep[Definition, page 3]{friedman1964partial}) and $\varGamma_{B,L_{b}}$
the fundamental solution of the backward parabolic equation $\left(\frac{\partial}{\partial t}+L_{b}\right)u=0$
(see \citep[Definition, page 27]{friedman1964partial}). 

1) The following holds: 
\begin{equation}
p_{b}(\tau,\xi,t,x)=\varGamma_{B,L_{b}}(\xi,\tau;x,t)=\varGamma_{F,L_{-b}-\nabla\cdot b}(x,t,\xi,\tau)\label{eq:fund-de-01}
\end{equation}
for all $0\leq\tau<t$ and $\xi,x\in\mathbb{R}^{d}$.

2) For given $\tau\geq0$, $\varphi$ and $g$, the function 
\begin{equation}
w(x,t)=\int p_{b}(\tau,\xi,t,x)\varphi(\xi)d\xi+\int_{\tau}^{t}\int p_{b}(\tau,\xi,t,x)g(\xi,s)d\xi ds\label{eq:for-w-sol1-2}
\end{equation}
solves the initial value problem of the parabolic equation: 
\begin{equation}
\left(\frac{\partial}{\partial t}-L_{-b}-\nabla\cdot b\right)w=g,\quad w(\cdot,0)=\varphi.\label{eq:for-adj1-2}
\end{equation}

3) If in addition $b(x,t)$ is solenoidal, i.e. $\nabla\cdot b=0$,
then 
\begin{equation}
p_{b}(\tau,\xi,t,x)=\varGamma_{F,L_{-b}}(x,t,\xi,\tau)\label{eq:fund-de-02}
\end{equation}
for all $0\leq\tau<t$ and $\xi,x\in\mathbb{R}^{d}$.
\end{lem}

The results in the previous lemma hold in fact under much weaker conditions
on $b$ and can be generalised to a large class of elliptic operators,
see \citep{aronson1967bounds,aronson1965uniqueness,friedman1964partial,stroock2008partial}
and other standard literature on parabolic equations for details. 

\subsection{An archetypical example}

Suppose the vorticity $\omega(x,t)$ of an incompressible fluid flow
with velocity $u(x,t)$, without applying external force, always lies
in the kernel of the rate-of-strain tensor, so that $\omega\cdot\nabla u=0$
identically, then the vorticity equation becomes
\begin{equation}
\left(\frac{\partial}{\partial t}-L_{-u(x,t)}\right)\omega(x,t)=0\label{eq:null-01}
\end{equation}
with the initial data $\omega(\cdot,0)=\omega_{0}$. Then, according
to Lemma \ref{lem:fund-sol-vortex}, 
\begin{equation}
\omega(x,t)=\int p_{u}(0,y,t,x)\omega_{0}(y)dy.\label{eq:null-02}
\end{equation}
On the other hand, since $\Delta u=-\nabla\wedge\omega$, according
to the Biot-Savart law, 
\begin{equation}
u(x,t)=\int G(x-z)\wedge\omega(z,t)dz,\label{eq:null-03}
\end{equation}
where $G(x)=-\frac{1}{4\pi}\frac{x}{|x|^{3}}$ is the vector valued
singular kernel in $\mathbb{R}^{3}$. Since $\omega$ is a solution
to \eqref{eq:null-01}, we are therefore able to rewrite the velocity
field \eqref{eq:null-03} in terms of the fundamental solution $p_{u}$,
to obtain that 
\begin{align}
u(x,t) & =\int\int G(x-z)\wedge\omega_{0}(y)p_{u}(0,y,t,z)dzdy\nonumber \\
 & =\int\left(\int G(x-z)p_{u}(0,y,t,z)dz\right)\wedge\omega_{0}(y)dy\nonumber \\
 & =\int\mathbb{E}\left[G(x-X(y,t))\right]\wedge\omega_{0}(y)dy,\label{eq:null-u}
\end{align}
where $X(y,t)$ is the Taylor diffusion process with infinitesimal
generator $L_{u}$ started at $y$ at $t=0$. That is, the solution
to the stochastic differential equation 
\begin{equation}
dX(y,t)=u(X(y,t),t)dt+\sqrt{2\nu}dB_{t},\quad X(y,0)=y,\label{eq:null-X}
\end{equation}
where $(B_{t})_{t\geq0}$ is a standard Brownian motion on a probability
space $(\Omega,\mathscr{F},\mathbb{P})$. Substituting \eqref{eq:null-u}
into \eqref{eq:null-X}, we may rewrite the previous stochastic differential
equation as 
\begin{equation}
dX(x,t)=\left.\left(\int\mathbb{E}\left[G\left(z-X(y,t)\right)\right]\wedge\omega_{0}(y)dy\right)\right|_{z=X(x,t)}dt+\sqrt{2\nu}dB,\quad X(x,0)=x,\label{eq:null-05-1}
\end{equation}
where $x$ runs through the state space $\mathbb{R}^{3}$. This is
the archetypical example of the SDEs we are going to study in the
present paper.

\subsection{Formulation of the problem}

Although our main examples come from the study of fluid dynamics,
it will be beneficial formulating the problem in a more general setting.
Still, we restrict our study to vector fields on Euclidean space $\mathbb{R}^{d}$.
Though the methods and the results can be generalised to tensor fields
with certain modifications.

Let $K(x)=(K_{j}^{i}(x))$ be a $d\times d$ matrix-valued `singular
integral' kernel, where $K_{j}^{i}$ are Borel measurable and locally
integrable. We are interested in the following stochastic differential
equation 
\begin{equation}
dX^{i}(x,t)=\left.\left(\int_{\mathbb{R}^{d}}\mathbb{E}\left[K_{j}^{i}\left(z-X(y,t)\right)\right]\omega_{0}^{j}(y)dy\right)\right|_{z=X(x,t)}dt+dB_{t}^{i},\quad X(x,0)=x,\label{eq:null-X-K}
\end{equation}
where $i=1,\ldots,d$, $\omega_{0}=(\omega_{0}^{i})$ is the initial
data, and $B=(B^{i})$ is a $d$-dimensional standard Brownian motion
on some probability space. Before we carry out a study of this class
of SDEs, let us reformulate \eqref{eq:null-X-K} in a different form
to facilitate our approach.

If $\mu$ is a measure on $(\mathbb{R}^{d},\mathcal{B}(\mathbb{R}^{d}))$,
then $K\star\mu=(K_{j}^{i}\star\mu)$ denotes the convolution of $K$
and the measure $\mu$ where
\begin{equation}
K_{j}^{i}\star\mu(x)=\int_{\R^{d}}K_{j}^{i}(x-y)\mu(dy)\label{eq:new2}
\end{equation}
for $i,j=1,\ldots,d$, as long as the right-hand side is well defined.

If $U$ is an $\mathbb{R}^{d}$-valued random variable on some probability
space $(\varOmega,\Fscr,\mathbb{P})$, then its distribution is denoted
by $\mathscr{L}(U)$. By definition 
\begin{equation}
K_{j}^{i}\star\mathscr{L}(U)(x)=\mathbb{E}\left[K_{j}^{i}(x-U)\right].\label{eq:new3}
\end{equation}
If, in addition, the law of $U$ has a pdf $p(x)$, then 
\begin{equation}
K_{j}^{i}\star\mathscr{L}(U)(x)=\int_{\mathbb{R}^{d}}K_{j}^{i}(x-y)p(y)dy\label{eq:new3-1}
\end{equation}
where the right-hand side is the convolution of $K$ and the function
$p$.

After having introduced the basic data $K$ and $\omega_{0}$ and
the notations, we are now in a position to reformulate the SDE we
are going to study:

\begin{equation}
dX^{i}(x,t)=\left[\int_{\R^{d}}\left(K_{j}^{i}\star\mathscr{L}(X(y,t))\right)(X(x,t))\omega_{0}^{j}(y)dy\right]dt+dB_{t}^{i}\label{eq:M_SDE1}
\end{equation}
with initial value $X(x,0)=x$ for $x\in\mathbb{R}^{d}$ and $i=1,\ldots,d$. 

The concepts of strong and weak solutions to \eqref{eq:M_SDE1} may
be defined accordingly. 

It will be convenient to introduce the following notations. If $Z=(Z(x,t))_{t\geq0}$
is a family of continuous processes on some probability space, which
is jointly continuous in $(x,t)$, then we may define a vector field
denoted by $b_{Z}$ whose components are given by
\begin{equation}
b_{Z}^{i}(x,t)=\int_{\R^{d}}\left(K_{j}^{i}\star\mathscr{L}(Z(y,t))\right)(x)\omega_{0}^{j}(y)dy,\label{eq:new6}
\end{equation}
$i=1,\ldots,d$. Notice that by definition, $b_{Z}$ depends only
on the one-dimensional marginal distributions of the process $(Z(y,t))_{t\geq0}$.

Suppose $b(x,t)$ is a time dependent vector field on $\mathbb{R}^{d}$,
we may define another $t$-dependent vector field on $\mathbb{R}^{d}$,
denoted by $K\diamond b(x,t)$, such that its $i$-th component is
given by
\begin{equation}
\int_{\R^{d}}\left(K_{j}^{i}\star\mathscr{L}(Z(y,t))\right)(x)\omega_{0}^{j}(y)dy,\label{eq:new5}
\end{equation}
where $Z=(Z(y,t))_{t\geq0}$ is the $L_{b}$-diffusion started at
$y$ at the moment $t=0$, so that $K\diamond b=b_{Z}$. Since $Z(y,t)$
has a transition probability density $p_{b}(0,y,t,z)$, we can write
\begin{equation}
(K\diamond b)^{i}(x,t)=\int_{\R^{d}}\left[\int_{\mathbb{R}^{d}}K_{j}^{i}(x-z)\omega_{0}^{j}(y)p_{b}(0,y,t,z)dz\right]dy.\label{eq:new5-1}
\end{equation}
We therefore define the mapping $K\diamond$ which sends a vector
field $b(x,t)$ to the vector field $K\diamond b(x,t)$. The non-linear
mapping $b\rightarrow K\diamond b$ will play a crucial r\^ole in
our study.

Under the above notations, we may rewrite SDE \eqref{eq:M_SDE1} as
\begin{equation}
dX(x,t)=b_{X}(X(x,t),t)dt+dB_{t},\quad X(x,0)=x,\label{eq:s-new7}
\end{equation}
for $x\in\mathbb{R}^{d}$. The following simple observation indeed
leads to the approach we are going to develop in what follows. 

\begin{lem}
\label{lem:2}Let $b(x,t)$ be a bounded and Borel measurable vector
field on $\mathbb{R}^{d}$, depending on time $t$. Suppose $K\diamond b=b$
and $(X,B)$ is a weak solution to the SDE 
\begin{equation}
dX^{i}(x,t)=b^{i}(X(x,t),t)dt+dB_{t}^{i},\quad X(x,0)=x,\label{eq:s-new9}
\end{equation}
where $B$ is a standard Brownian motion on a probability space and
$i=1,\ldots,d$. Then 
\begin{equation}
dX^{i}(x,t)=\left[\int_{\R^{d}}\left(K_{j}^{i}\star\mathscr{L}(X(y,t))\right)(X(x,t))\omega_{0}^{j}(y)dy\right]dt+dB_{t}^{i}\label{eq:s new10}
\end{equation}
where $i=1,\ldots,d$. That is, $(X,B)$ is a weak solution to \eqref{eq:M_SDE1}. 
\end{lem}

This lemma follows by definition: since $K\diamond b=b$, so $b=b_{X}$,
which implies that 
\[
\int_{\R^{d}}\left(K_{j}^{i}\star\mathscr{L}(X(y,t))\right)(x)\omega_{0}^{j}(y)dy=b^{i}(x,t)
\]
and therefore \eqref{eq:s new10} follows from \eqref{eq:s-new9}
immediately.
\begin{example}
If $d=3$ and $K=(K^{1},K^{2},K^{3})$, then we set $K_{j}^{i}=\varepsilon^{ikj}K^{k}$
and SDE \eqref{eq:s new10} becomes
\[
dX(x,t)=\left.\left(\int_{\mathbb{R}^{d}}\mathbb{E}\left[K\left(z-X(y,t)\right)\right]\wedge\omega_{0}(y)dy\right)\right|_{z=X(x,t)}dt+dB_{t},\quad X(x,0)=x,
\]
which is the random vortex dynamical model, where $\omega_{0}$ represents
the initial vorticity.
\end{example}

\section{Several facts about diffusions with bounded drifts}

In this section, we collect a few facts on diffusion processes and
prove several technical estimates which will be used later in next
section. Let $b(x,t)$ be a Borel measurable vector field on the Euclidean
space $\mathbb{R}^{d}$, dependent on the time parameter $t\geq0$.
It is assumed that $b(x,t)$ is bounded: $|b(x,t)|\leq A$ for every
$x$ and $t$, where $A$ is a non-negative constant. Then the unique
$L_{b}$-diffusion (in the sense of weak solutions) may be constructed
by using Cameron-Martin formula (see \citep[Theorem 6.4.2, page 154]{stroock1979multidimensional}).
Let $B=(B_{t})_{t\geq0}$ be a $d$-dimensional standard Brownian
motion on $(\varOmega,\mathcal{F},\mathbb{P})$. Let $\mathcal{F}_{t}=\sigma\left\{ B_{s}:s\leq t\right\} $
and $\mathcal{F}_{\infty}=\sigma\left\{ B_{s}:s\geq0\right\} $ be
the natural filtration generated by this Brownian motion. Given $x$
and $\tau\geq0$, define the exponential martingale called the Cameron-Martin
density
\begin{equation}
R_{b}(\tau,x,t)=e^{N_{b}(\tau,x,t)}\label{eq:s-new1}
\end{equation}
where, for simplicity, we have written
\begin{equation}
N_{b}(\tau,x,t):=\int_{\tau}^{t}b(r,B_{r}-B_{\tau}+x)dB_{r}-\frac{1}{2}\int_{\tau}^{t}|b|^{2}(r,B_{r}-B_{\tau}+x)dr\label{eq:s-new3}
\end{equation}
for $t\geq\tau$. If $\tau=0$, then the symbol $\tau$ will be suppressed
from the notations. Next, construct the probability $\mathbb{P}^{\tau,x}$
on $(\varOmega,\mathcal{F}_{\infty})$ such that
\begin{equation}
\left.\frac{d\mathbb{P}^{\tau,x}}{d\mathbb{P}}\right|_{\mathcal{F}_{t}}=R_{b}(\tau,x,t)\label{eq:2-new1}
\end{equation}
for all $t\geq\tau$. Then the family $\{\mathbb{P}^{\tau,x}:\tau\geq0,x\in\mathbb{R}^{d}\}$
on $(\varOmega,\mathcal{F}_{\infty})$ is a diffusion family with
generator $L_{b}$ (see for example \citep{Karatzas and Shreve 1988,stroock1979multidimensional}).
In particular
\begin{equation}
\int_{\mathbb{R}^{d}}f(y)p_{b}(\tau,x,t,y)dy=\mathbb{P}\left[R_{b}(\tau,x,t)f(B_{t}-B_{\tau})\right]\label{eq:2-new2}
\end{equation}
for any Borel function $f$, as long as one of the integrals in the
equation makes sense, where\emph{ $p_{b}(\tau,x,t,y)$} is the transition
probability density function of the $L_{b}$-diffusion. It is known
(see for example \citep{stroock2008partial}) that $p_{b}(\tau,x,t,y)$
is positive and continuous on any $t>\tau\geq0$ and $x,y\in\mathbb{R}^{d}$.
Moreover, for every $T>0$, there is a constant $M$ depending on
$A$, $d$ and $T$ only, such that 
\begin{equation}
\frac{1}{Mt^{d/2}}e^{-M\frac{|y-x|^{2}}{t}}\leq p_{b}(\tau,x,\tau+t,y)\leq\frac{M}{t^{d/2}}e^{-\frac{|y-x|^{2}}{Mt}}\label{eq:aronson est1}
\end{equation}
for all $\tau\geq0$ and $T\geq t>0$. This is the so-called Aronson
estimate (see \citep{Aronson 1968-main,stroock1988diffusion,stroock2008partial}
for example). In our study, we need more precise information about
the constant $M$, which was obtained in \citep{qian2002sharp,qian2003comparison}.
\begin{lem}
\label{lem:qian-sharp-bd} There is a positive universal constant
$\kappa$, depending only on the dimension $d$ and $1<q<\frac{d}{d-1}$,
such that
\begin{equation}
p_{b}(\tau,x,\tau+t,y)\leq\frac{1}{(2\pi t)^{d/2}}e^{-\frac{|x-y|^{2}}{2t}}\left(1+\kappa A\left(\sqrt{t}+|x-y|\right)e^{\frac{q-1}{2qt}|x-y|^{2}+\frac{A^{2}t}{2(q-1)}}\right)\label{eq:new1}
\end{equation}
for all $x,y\in\R^{d}$, $\tau\geq0$ and $t>0$.
\end{lem}

As a consequence, we establish the following estimate, which will
play a crucial r\^ole in the proof of our main theorem.
\begin{lem}
\label{lem:cut-small}Let $f\in L^{\infty}(\mathbb{R}^{d})$, $\gamma\in[0,d)$
and $\rho>0$. Define
\begin{equation}
I(f,x,t,\rho,\gamma)=\int_{\mathbb{R}^{d}}\int_{\{|z|<\rho\}}\frac{1}{|z|^{\gamma}}|f(y)|p_{b}(0,y,t,x-z)dzdy\label{eq:cu-est1}
\end{equation}
for any $x\in\mathbb{R}^{d}$ and $t>0$. Then there exists a universal
positive constant $\kappa_{1}$ depending only on $d$, such that
\begin{equation}
I(f,x,t,\rho,\gamma)\leq\frac{\rho^{d-\gamma}}{d-\gamma}\kappa_{1}\left\Vert f\right\Vert _{\infty}\left(1+A\sqrt{t}e^{\frac{A^{2}t}{2(q-1)}}\right)\label{eq:key-1est}
\end{equation}
for all $x$ and $t>0$.
\end{lem}

\begin{proof}
Without losing generality, we may assume that $f\geq0$. Using the
sharp estimate \eqref{eq:new1}, we have
\begin{align*}
I(f,x,t,\rho,\gamma) & \leq\int_{\mathbb{R}^{d}}\int_{\{|z|<\rho\}}\frac{f(y)}{|z|^{\gamma}}\frac{e^{-\frac{|y-x+z|^{2}}{2t}}}{(2\pi t)^{d/2}}\\
 & \quad\times\left(1+\kappa A\left(\sqrt{t}+|y-x+z|\right)e^{\frac{q-1}{2qt}|y-x+z|^{2}+\frac{A^{2}t}{2(q-1)}}\right)dzdy\\
 & =\int_{\mathbb{R}^{d}}\int_{\{|z|<\rho\}}\frac{f(y+x-z)}{|z|^{\gamma}}\frac{e^{-\frac{|y|^{2}}{2t}}}{(2\pi t)^{d/2}}\left(1+\kappa A\left(\sqrt{t}+|y|\right)e^{\frac{q-1}{2qt}|y|^{2}+\frac{A^{2}t}{2(q-1)}}\right)dzdy\\
 & \leq\left\Vert f\right\Vert _{\infty}\left(\int_{\{|z|<\rho\}}\frac{1}{|z|^{\gamma}}dz\right)\left(1+\kappa A\sqrt{t}e^{\frac{A^{2}t}{2(q-1)}}\int_{\mathbb{R}^{d}}\frac{e^{-\frac{|y|^{2}}{2t}}}{(2\pi t)^{d/2}}\left(1+\frac{|y|}{\sqrt{t}}\right)e^{\frac{q-1}{2qt}|y|^{2}}dy\right)\\
 & =\left\Vert f\right\Vert _{\infty}\left(\int_{\{|z|<\rho\}}\frac{1}{|z|^{\gamma}}dz\right)\left(1+\kappa A\sqrt{t}e^{\frac{A^{2}t}{2(q-1)}}\int_{\mathbb{R}^{d}}\frac{e^{-\frac{|y|^{2}}{2}}}{(2\pi)^{d/2}}\left(1+|y|\right)e^{\frac{q-1}{2q}|y|^{2}}dy\right)\\
 & \leq\frac{\rho^{d-\gamma}}{d-\gamma}\kappa_{1}\left\Vert f\right\Vert _{\infty}\left(1+A\sqrt{t}e^{\frac{A^{2}t}{2(q-1)}}\right)
\end{align*}
where
\[
\kappa_{1}=\max\left\{ \textrm{Vol}(S^{d-1}),\kappa\int_{\mathbb{R}^{d}}\left(1+|y|\right)\frac{e^{-\frac{|y|^{2}}{2q}}}{(2\pi)^{d/2}}dy\right\} 
\]
and the proof is complete.
\end{proof}
We also need the following estimate which is completely elementary
though.
\begin{lem}
\label{lem:cut-big}Let $f\in L^{1}(\mathbb{R}^{d})$ and $\rho>0$,
$\gamma\geq0$ be two constants. Define
\begin{equation}
J(f,x,t,\rho,\gamma)=\int_{\mathbb{R}^{d}}\int_{\{|z|\geq\rho\}}\frac{1}{|z|^{\gamma}}|f(y)|p_{b}(0,y,t,x-z)dzdy\label{eq:def-J1}
\end{equation}
for all $x\in\mathbb{R}^{d}$ and $t>0$. Then
\begin{equation}
J(f,x,t,\rho,\gamma)\leq\frac{1}{\rho^{\gamma}}\left\Vert f\right\Vert _{L^{1}}\label{eq:J-out eq1}
\end{equation}
for all $x$ and $t>0$, $\rho>0$ and $\gamma\geq0$.
\end{lem}

\begin{proof}
Since $\gamma\geq0$, so that
\begin{align*}
J(f,x,t,\rho,\gamma) & \leq\frac{1}{\rho^{\gamma}}\int_{\mathbb{R}^{d}}\int_{\{|z|\geq\rho\}}|f(y)|p_{b}(0,y,t,x-z)dzdy\\
 & \leq\frac{1}{\rho^{\gamma}}\int_{\mathbb{R}^{d}}\int_{\mathbb{R}^{d}}|f(y)|p_{b}(0,y,t,x-z)dzdy\\
 & =\frac{1}{\rho^{\gamma}}\int_{\mathbb{R}^{d}}\int_{\mathbb{R}^{d}}|f(y)|p_{b}(0,y,t,z)dzdy\\
 & =\frac{1}{\rho^{\gamma}}\int_{\mathbb{R}^{d}}|f(y)|dy\\
 & =\frac{1}{\rho^{\gamma}}\left\Vert f\right\Vert _{L^{1}}
\end{align*}
and the proof is complete.
\end{proof}

\section{Weak solutions}

In this section we prove under certain conditions that there is a
unique weak solution to \eqref{eq:M_SDE1}. To this end we make several
assumptions on $\omega_{0}$ and $K$ which will be in force throughout
the remainder of the paper.

Let $C_{0},C_{1}$ and $C_{\infty}$ be three non-negative constants.
It is assumed that $K=(K_{j}^{i})$ satisfies the following growth
condition: there are two constants $\gamma_{1}\in[0,d)$ and $\gamma_{2}\geq0$
such that
\begin{equation}
\left|K(x)\right|\leq\frac{C_{0}}{|x|^{\gamma_{1}}}\quad\textrm{ for all }x\neq0\textrm{ and }|x|<1,\label{eq:new7}
\end{equation}
and
\begin{equation}
\left|K(x)\right|\leq\frac{C_{0}}{|x|^{\gamma_{2}}}\quad\textrm{ for all }|x|\geq1.\label{eq:new7-b}
\end{equation}

In addition we assume that the initial vorticity $\omega_{0}$ is
bounded and integrable such that $\left\Vert \omega_{0}\right\Vert _{L^{1}}\leq C_{1}$
and $\left\Vert \omega_{0}\right\Vert _{\infty}\leq C_{\infty}$.
Choose and fix a number $q\in(1,\frac{d}{d-1})$ and set
\begin{equation}
C_{K}=C_{0}\left(\frac{\kappa_{1}C_{\infty}}{d-\gamma_{1}}\left(1+e^{\frac{1}{2(q-1)}}\right)+C_{1}\right)\textrm{ and }\quad T_{K}=\frac{1}{C_{K}^{2}}.\label{eq:new 8}
\end{equation}
The crucial fact about $C_{K}$ and $T_{K}$ is that they depend on
$C_{0},C_{1}$, $C_{\infty}$ and $\gamma_{1}$ only.
\begin{lem}
\label{lem:b-unif-bdd} If $b(x,t)$ is a time-dependent vector field
 such that $|b(x,t)|\leq C_{K}$ for all $x\in\R^{d}$ and
$t\leq T_{K}$, then $K\diamond b$ is also bounded with the same
bound. That is $|K\diamond b(x,t)|\leq C_{K}$ for all $x\in\R^{d}$
and $t\leq T_{K}$. 
\end{lem}

\begin{proof}
Let $B=\left\{ z\in\mathbb{R}^{d}:|z|<1\right\} $. Then 
\begin{align}
\left|K\diamond b(x,t)\right| & =\bigg|\int_{\mathbb{R}^{d}}\int_{B}K(z)\omega_{0}(y)p_{b}(0,y,t,x-z)dzdy\nonumber \\
 & \quad+\int_{\mathbb{R}^{d}}\int_{\R^{d}\backslash B}K(z)\omega_{0}(y)p_{b}(0,y,t,x-z)dzdy\bigg|\nonumber \\
 & \leq\int_{\mathbb{R}^{d}}\int_{B}\left|K(z)\right|\left|\omega_{0}(y)\right|p_{b}(0,y,t,x-z)dzdy\nonumber \\
 & \quad+\int_{\mathbb{R}^{d}}\left|\omega_{0}(y)\right|\int_{\R^{d}\backslash B}|K(z)|p_{b}(0,y,t,x-z)dzdy\nonumber \\
 & =:I_{1}+I_{2}\label{eq:b-I1-I2}
\end{align}
The estimate for $I_{1}$ follows directly from Lemma \ref{lem:cut-small}:
\begin{align*}
I_{1} & \leq C_{0}I(\omega_{0},x,t,1,\gamma_{1})\\
 & \leq\frac{\kappa_{1}C_{0}C_{\infty}}{d-\gamma_{1}}\left(1+A\sqrt{t}e^{\frac{A^{2}t}{2(q-1)}}\right)
\end{align*}
where $A=C_{K}$. By Lemma \ref{lem:cut-big} we deduce that
\[
I_{2}\leq C_{0}J(\omega_{0},x,t,1,\gamma_{2})\leq C_{0}C_{1}.
\]
Putting the estimates for $I_{1}$ and $I_{2}$ together we may conclude
that 
\[
\left|K\diamond b(x,t)\right|\leq\frac{\kappa_{1}C_{0}C_{\infty}}{d-\gamma_{1}}\left(1+A\sqrt{t}e^{\frac{A^{2}t}{2(q-1)}}\right)+C_{0}C_{1}
\]
for all $x\in\R^{d}$ and $t\geq0$. Since $A\sqrt{t}\leq1$ for any
$t\leq T_{K}$, we therefore have
\[
\left|K\diamond b(x,t)\right|\leq\frac{\kappa_{1}C_{0}C_{\infty}}{d-\gamma_{1}}\left(1+e^{\frac{1}{2(q-1)}}\right)+C_{0}C_{1}
\]
for any $x$ and $t\leq T_{K}$. The conclusion then follows immediately
from the definition of $C_{K}$ and $T_{K}$.
\end{proof}
Next we are going to establish another key estimate for the mapping
$b\rightarrow K\diamond b$, where $b(x,t)$ are vector fields such
that $|b(x,t)|\leq C_{K}$ for any $t\leq T_{K}$.
\begin{lem}
\label{lem:6} There exists a positive constant $C_{L}$ depending
only \textup{on $C_{0},C_{1},C_{\infty}$,} such that for any $b(x,t)$
and $\tilde{b}(x,t)$ satisfying that $|b(x,t)|\leq C_{K}$ and $|\tilde{b}(x,t)|\leq C_{K}$
for all $x$ and $t\leq T_{K}$ we have
\begin{equation}
|K\diamond b(x,t)-K\diamond\tilde{b}(x,t)|\leq\left(t+\sqrt{t}\right)C_{L}\left\Vert b-\tilde{b}\right\Vert _{L^{\infty}(\mathbb{R}^{d},[0,t])}\label{eq:key-2}
\end{equation}
for all $x$ and $t\leq T_{K}$. 
\end{lem}

\begin{proof}
We prove this by using Cameron-Martin formula \citep[Theorem 6.4.2, page 154]{stroock1979multidimensional}.
Let $B$ be a $d$-dimensional standard Brownian motion on some probability
space $(\varOmega,\mathcal{F},\mathbb{P})$ and $R_{c}(x,t)=e^{N_{c}(x,t)}$
be the Cameron-Martin density (see \eqref{eq:s-new1} and \eqref{eq:s-new3})
with respect to the vector field $c$ starting at $x$ at the moment
$0$. Then
\begin{align*}
N_{b}(x,t)-N_{\tilde{b}}(x,t) & =\int_{0}^{t}\left(b-\tilde{b}\right)(r,B_{r}+x)dB_{r}\\
 & \quad-\frac{1}{2}\int_{0}^{t}\left(|b|^{2}-|\tilde{b}|^{2}\right)(r,B_{r}+x)dr\\
 & =M(x,t)+A(x,t),
\end{align*}
where 
\[
M(x,t):=\int_{0}^{t}\left(b-\tilde{b}\right)(r,B_{r}+x)dB_{r},
\]
whose quadratic process 
\begin{equation}
\left\langle M\right\rangle _{t}=\int_{0}^{t}|b-\tilde{b}|^{2}(r,B_{r}+x)dr\leq\left\Vert b-\tilde{b}\right\Vert _{\infty}^{2}t,\label{eq:b-new1}
\end{equation}
and 
\begin{align*}
A(x,t) & :=-\frac{1}{2}\int_{0}^{t}\left(|b|^{2}-|\tilde{b}|^{2}\right)(r,B_{r}+x)dr\\
 & =-\frac{1}{2}\int_{0}^{t}\left(b-\tilde{b}\right)\left(b+\tilde{b}\right)(r,B_{r}+x)dr.
\end{align*}
It is clear that 
\begin{equation}
|A(x,t)|\leq C_{K}t\left\Vert b-\tilde{b}\right\Vert _{L^{\infty}([0,t])}\label{eq:s-new11}
\end{equation}
and therefore,
\begin{equation}
\left|N_{b}(x,t)-N_{\tilde{b}}(x,t)\right|\leq|M_{t}|+C_{K}t\left\Vert b-\tilde{b}\right\Vert _{L^{\infty}([0,t])}.\label{eq:s new 12}
\end{equation}
Now we write
\begin{align}
R_{b}-R_{\tilde{b}} & =\int_{0}^{1}\frac{d}{ds}e^{(1-s)N_{\tilde{b}}+sN_{b}}ds\nonumber \\
 & =\int_{0}^{1}e^{(1-s)N_{\tilde{b}}+sN_{b}}(N_{b}-N_{\tilde{b}})ds.\label{eq:ad-01}
\end{align}
Let $b_{s}=(1-s)\tilde{b}+sb$ for $s\in[0,1]$. Then $|b_{s}(x,t)|\leq C_{K}$
and
\[
(1-s)N_{\tilde{b}}+sN_{b}=N_{b_{s}}-s(1-s)|b-\tilde{b}|^{2}.
\]
Substituting this equality into \eqref{eq:ad-01} to obtain
\begin{equation}
R_{b}-R_{\tilde{b}}=\int_{0}^{1}R_{b_{s}}e^{-(1-s)s|b-\tilde{b}|^{2}}(N_{b}-N_{\tilde{b}})ds.\label{eq:s-new2}
\end{equation}
Now we are in a position to study the non-linear mapping $c\rightarrow K\diamond c$.
According to Cameron-Martin formula
\begin{align*}
K\diamond c(x,t) & =\int_{\R^{d}}\left[\int_{\mathbb{R}^{d}}K(x-z)\omega_{0}(y)p_{c}(0,y,t,z)dz\right]dy\\
 & =\int_{\R^{d}}\mathbb{E}\left[R_{c}(y,t)K(x-B_{t}-y)\right]\omega_{0}(y)dy.
\end{align*}
Split $K$ into a sum $K=K_{1}+K_{2}$ where
\[
K_{1}(z)=1_{\left\{ |z|<1\right\} }K(z)\textrm{ and }K_{2}(z)=1_{\left\{ |z|\geq1\right\} }K(z).
\]
Then
\begin{align*}
K\diamond c(x,t) & =\int_{\R^{d}}\mathbb{E}\left[R_{c}(y,t)K_{1}(x-B_{t}-y)\right]\omega_{0}(y)dy\\
 & \quad+\int_{\R^{d}}\mathbb{E}\left[R_{c}(y,t)K_{2}(x-B_{t}-y)\right]\omega_{0}(y)dy.
\end{align*}
Let 
\[
D(x,t):=K\diamond b(x,t)-K\diamond\tilde{b}(x,t).
\]
Then by using the previous formula for $K\diamond c$ we have
\begin{align*}
D(x,t) & =\int_{\R^{d}}\mathbb{E}\left[\left(R_{b}(y,t)-R_{\tilde{b}}(y,t)\right)K_{1}(x-B_{t}-y)\right]\omega_{0}(y)dy\\
 & \quad+\int_{\R^{d}}\mathbb{E}\left[\left(R_{b}(y,t)-R_{\tilde{b}}(y,t)\right)K_{2}(x-B_{t}-y)\right]\omega_{0}(y)dy\\
 & =:J_{1}+J_{2}.
\end{align*}
Substituting \eqref{eq:s-new2} into $J_{1}$, we may write
\[
J_{1}=\int_{\R^{d}}\mathbb{E}\left[\left(\int_{0}^{1}R_{b_{s}}(y,t)e^{-(1-s)s|b-\tilde{b}|^{2}}(N_{b}(y,t)-N_{\tilde{b}}(y,t))ds\right)K_{1}(x-B_{t}-y)\right]\omega_{0}(y)dy
\]
and therefore
\begin{align*}
\left|J_{1}\right| & \leq\int_{0}^{1}\int_{\R^{d}}\mathbb{E}\left[R_{b_{s}}(y,t)\left|N_{b}(y,t)-N_{\tilde{b}}(y,t)\right|\left|K_{1}(x-B_{t}-y)\right|\right]\left|\omega_{0}(y)\right|dyds\\
 & \leq C_{0}\int_{0}^{1}\int_{\R^{d}}\mathbb{E}\left[R_{b_{s}}(y,t)\frac{1_{\left\{ |x-B_{t}-y|<\rho\right\} }}{|x-B_{t}-y|^{\gamma_{1}}}|M_{t}|\right]\left|\omega_{0}(y)\right|dyds\\
 & \quad+C_{0}C_{K}t\left\Vert b-\tilde{b}\right\Vert _{L^{\infty}([0,t])}\int_{0}^{1}\int_{\R^{d}}\mathbb{E}\left[R_{b_{s}}(y,t)\frac{1_{\left\{ |x-B_{t}-y|<\rho\right\} }}{|x-B_{t}-y|^{\gamma_{1}}}\right]\left|\omega_{0}(y)\right|dyds\\
 & =:J_{1,1}+J_{1,2},
\end{align*}
where the second inequality comes from \eqref{eq:s new 12}.

To deal with $J_{1,1}$, choose and fix $\alpha,\beta>1$ such that
$\alpha\gamma_{1}<d$ and $\alpha^{-1}+\beta^{-1}=1$. Then
\[
R_{b_{s}}^{\alpha}=R_{\alpha b_{s}}e^{\frac{\alpha^{2}-1}{2}\int_{0}^{t}|b_{s}|^{2}dr}\leq e^{\frac{\alpha^{2}-1}{2}C_{K}^{2}t}R_{\alpha b_{s}}.
\]
Also, by using Burkholder-Davis-Gundy inequality (see for example
\citep[Theorem 3.1, page 110]{Ikeda Watanabe}),
\begin{equation}
\mathbb{E}\left[|M_{t}|^{\beta}\right]\leq C_{\beta}\mathbb{E}\left[\left\langle M\right\rangle _{t}^{\frac{\beta}{2}}\right]\leq C_{\beta}t^{\frac{\beta}{2}}\left\Vert b-\tilde{b}\right\Vert _{L^{\infty}(\mathbb{R}^{d}\times[0,t])}^{\beta}.\label{eq:M-quadratic1}
\end{equation}
Thus, by applying H\"older's inequality in $J_{1,1}$, we deduce that
\begin{align*}
J_{1,1} & =C_{0}\int_{0}^{1}\int_{\R^{d}}\mathbb{E}\left[R_{b_{s}}(y,t)\frac{1_{\left\{ |x-B_{t}-y|<1\right\} }}{|x-B_{t}-y|^{\gamma_{1}}}|M_{t}|\right]\left|\omega_{0}(y)\right|dyds\\
 & \leq C_{0}\int_{0}^{1}\sqrt[\alpha]{\int_{\R^{d}}\mathbb{E}\left[R_{b_{s}}^{\alpha}(y,t)\frac{1_{\left\{ |x-B_{t}-y|<1\right\} }}{|x-B_{t}-y|^{\alpha\gamma_{1}}}\right]\left|\omega_{0}(y)\right|dy}\\
 & \quad\times\sqrt[\beta]{\int_{\R^{d}}\mathbb{E}\left[|M_{t}|^{\beta}\right]\left|\omega_{0}(y)\right|dy}ds\\
 & \leq C_{0}\sqrt[\beta]{C_{1}C_{\beta}}e^{\frac{\alpha^{2}-1}{2\alpha}C_{K}^{2}t}\sqrt{t}\left\Vert b-\tilde{b}\right\Vert _{L^{\infty}([0,t])}\\
 & \quad\times\int_{0}^{1}\sqrt[\alpha]{\int_{\R^{d}}\mathbb{E}\left[R_{\alpha b_{s}}(y,t)\frac{1_{\left\{ |x-B_{t}-y|<1\right\} }}{|x-B_{t}-y|^{\alpha\gamma_{1}}}\right]\left|\omega_{0}(y)\right|dy}ds\\
 & =C_{0}\sqrt[\beta]{C_{1}C_{\beta}}e^{\frac{\alpha^{2}-1}{2\alpha}C_{K}^{2}t}\sqrt{t}\left\Vert b-\tilde{b}\right\Vert _{L^{\infty}([0,t])}\\
 & \quad\times\int_{0}^{1}\sqrt[\alpha]{\int_{\R^{d}}\int_{|z|<1}\frac{1}{|z|^{\alpha\gamma_{1}}}p_{\alpha b_{s}}(0,y,t,x-z)dz\left|\omega_{0}(y)\right|dy}ds\\
 & \leq C_{0}\sqrt[\beta]{C_{1}C_{\beta}}e^{\frac{\alpha^{2}-1}{2\alpha}C_{K}^{2}t}\sqrt[\alpha]{\frac{\kappa_{1}\left(1+\alpha C_{K}\sqrt{t}e^{\frac{\alpha^{2}C_{K}^{2}t}{2(q-1)}}\right)}{d-\alpha\gamma_{1}}}\sqrt{t}\left\Vert b-\tilde{b}\right\Vert _{L^{\infty}(\mathbb{R}^{d}\times[0,t])}\\
 & \leq C_{0}\sqrt[\beta]{C_{1}C_{\beta}}e^{\frac{\alpha^{2}-1}{2\alpha}}\sqrt[\alpha]{\frac{\kappa_{1}\left(1+\alpha e^{\frac{\alpha^{2}}{2(q-1)}}\right)}{d-\alpha\gamma_{1}}}\sqrt{t}\left\Vert b-\tilde{b}\right\Vert _{L^{\infty}(\mathbb{R}^{d}\times[0,t])},
\end{align*}
 for all $t\leq T_{K}$, where the second inequality follows from
\eqref{eq:M-quadratic1}, the third inequality follows from Lemma
\ref{lem:cut-small}, and the last inequality follows from the fact
that $C_{K}t\leq1$ for all $t\leq T_{K}$.

To deal with $J_{1,2}$, we apply Lemma \ref{lem:cut-big} and obtain
that
\begin{align*}
J_{1,2} & =C_{0}C_{K}t\left\Vert b-\tilde{b}\right\Vert _{L^{\infty}([0,t])}\\
 & \quad \times\int_{0}^{1}\int_{\R^{d}}\mathbb{E}\left[R_{b_{s}}(y,t)\frac{1_{\left\{ |x-B_{t}-y|<1\right\} }}{|x-B_{t}-y|^{\gamma_{1}}}\right]\left|\omega_{0}(y)\right|dyds\\
 & =C_{0}C_{K}t\left\Vert b-\tilde{b}\right\Vert _{L^{\infty}([0,t])}\\
 & \quad  \times\int_{0}^{1}\int_{\R^{d}}\left(\int_{\{|z|<1\}}\frac{1}{|z|^{\gamma_{1}}}p_{b_{s}}(0,y,t,x-z)dz\right)\left|\omega_{0}(y)\right|dyds\\
 & \leq\frac{C_{0}C_{\infty}C_{K}\kappa_{1}}{d-\gamma_{1}}\left(1+C_{K}\sqrt{t}e^{\frac{C_{K}^{2}t}{2(q-1)}}\right)t\left\Vert b-\tilde{b}\right\Vert _{L^{\infty}(\mathbb{R}^{d}\times[0,t])}\\
 & \leq\frac{C_{0}C_{\infty}C_{K}\kappa_{1}}{d-\gamma_{1}}\left(1+e^{\frac{1}{2(q-1)}}\right)t\left\Vert b-\tilde{b}\right\Vert _{L^{\infty}(\mathbb{R}^{d}\times[0,t])}
\end{align*}
for any $t\leq T_{K}$, where the first inequality follows from the
estimate in Lemma \ref{lem:cut-big}.

Now we treat with $J_{2}$. Since
\[
J_{2}=\int_{\R^{d}}\mathbb{E}\left[\left(\int_{0}^{1}R_{b_{s}}(y,t)e^{-(1-s)s|b-\tilde{b}|^{2}}(N_{b}(y,t)-N_{\tilde{b}}(y,t))ds\right)K_{2}(x-B_{t}-y)\right]\omega_{0}(y)dy
\]
and $|K_{2}(z)|\leq C_{0}$, so by \eqref{eq:s new 12} we have
\begin{align*}
|J_{2}| & \leq C_{0}\int_{0}^{1}\int_{\R^{d}}\mathbb{E}\left[R_{b_{s}}(y,t)|M_{t}|\right]\left|\omega_{0}(y)\right|dyds\\
 & \quad +C_{0}C_{1}C_{K}t\left\Vert b-\tilde{b}\right\Vert _{L^{\infty}(\mathbb{R}^{d}\times[0,t])}\\
 & \leq C_{0}\int_{0}^{1}\int_{\R^{d}}\sqrt{\mathbb{E}\left[|M_{t}|^{2}\right]}\sqrt{\mathbb{E}\left[R_{b_{s}}^{2}(y,t)\right]}\left|\omega_{0}(y)\right|dyds\\
 & \quad +C_{0}C_{1}C_{K}t\left\Vert b-\tilde{b}\right\Vert _{L^{\infty}(\mathbb{R}^{d}\times[0,t])}\\
 & \leq C_{0}C_{1}\sqrt{e^{\frac{\alpha^{2}-1}{2}}}\sqrt{t}\left\Vert b-\tilde{b}\right\Vert _{L^{\infty}([0,t])}\\
 & \quad +C_{0}C_{1}C_{K}t\left\Vert b-\tilde{b}\right\Vert _{L^{\infty}(\mathbb{R}^{d}\times[0,t])},
\end{align*}
where the last inequality comes from \eqref{eq:M-quadratic1}. Putting
these estimates for $J_{1}$ and $J_{2}$ together, we deduce \eqref{eq:key-2}
with a positive constant $C_{L}$ which depends only on the structure
constants $C_{0}$, $C_{1}$, $C_{\infty}$,$\gamma_{1}$ and $d$
(as $\alpha$ and $q$ are constants depending only on $d$ and $\gamma_{1}$).
For example
\begin{align}
C_{L} & =C_{0}C_{1}\sqrt{e^{\frac{\alpha^{2}-1}{2}}}(C_{K}+1)\nonumber \\
 & \quad +C_{0}C_{\infty}C_{K}\frac{1}{d-\gamma_{1}}\kappa_{1}\left(1+e^{\frac{1}{2(q-1)}}\right)\\
 & \quad +C_{0}\sqrt[\beta]{C_{1}C_{\beta}}e^{\frac{\alpha^{2}-1}{2\alpha}}\sqrt[\alpha]{\frac{1}{d-\alpha\gamma_{1}}\kappa_{1}\left(1+\alpha e^{\frac{\alpha^{2}}{2(q-1)}}\right)}\label{eq:C_L}
\end{align}
will do. 
\end{proof}
We are now in a position to prove the main result about weak solutions
to \eqref{eq:M_SDE1}.
\begin{thm}
\label{thm:main 1} There exist two positive constants $T_{L}$ and
$C_{K}$ depending on $C_{0},C_{1},C_{\infty}$, $\gamma_{1}\in[0,d)$
and $d$ only such that the followings hold:

1) The (non-linear) mapping $b\rightarrow K\diamond b$ is contractive
on the space of bounded time-dependent vector fields. More precisely
\begin{equation}
\left\Vert K\diamond b-K\diamond\tilde{b}\right\Vert _{\infty}\leq\frac{1}{2}\left\Vert b-\tilde{b}\right\Vert _{\infty}\label{eq:contr1}
\end{equation}
for any vector fields $b$ and $\tilde{b}$ such that $\left\Vert b\right\Vert _{\infty}\leq C_{K}$,
$\Vert\tilde{b}\Vert_{\infty}\leq C_{K}$, where $\left\Vert b\right\Vert _{\infty}=\left\Vert b\right\Vert _{L^{\infty}(\mathbb{R}^{d}\times[0,T_{L}])}$.
Hence, there is a unique $b$ such that $K\diamond b=b$. 

2) There is a unique weak solution $(X,B)$ on some probability space
to the SDE \eqref{eq:M_SDE1} up to time $T_{L}$, where $B$ is a
Brownian motion and $X$ satisfies \eqref{eq:M_SDE1}, and 
\[
b^{i}(\cdot,t)=\int_{\R^{d}}\left(K_{j}^{i}\star\mathscr{L}(X(y,t))\right)\omega_{0}^{j}(y)dy
\]
is bounded for $i=1,\cdots,d$.
\end{thm}

\begin{proof}
Choose $T_{L}=\frac{1}{4}C_{L}\wedge1$ where $C_{L}$ is given by
\eqref{eq:C_L}. Then \eqref{eq:contr1} follows immediately. 

The second part then follows from Lemma \ref{lem:2}.
\end{proof}
We finish this section with a comment on the global solutions of \eqref{eq:M_SDE1}.
As long as $K$ is a singular integral kernel, bounded at infinity,
we have shown that there is a unique weak solution to \eqref{eq:M_SDE1}
with initial data $\omega_{0}\in L^{1}(\mathbb{R}^{d})\cap L^{\infty}(\mathbb{R}^{d})$
for the time duration $[0,T_{L}]$, where $T_{L}$ depends only on
the structure constants $C_{i}$ ($i=0,1,\infty$) and $\gamma_{1}\in[0,d)$
. However, we are unable to conclude that the weak solution exists
for all time $t$. The reason is that the SDE \eqref{eq:M_SDE1} does
not define a dynamical system, which is not proposed as an initial
value problem. 

Finally, we should point out that we do not claim, although we strongly
believe it is not the case, if $\omega_{0}$ and $K$ are regular
enough, there are other fixed vector fields $c$ in the sense that
$K\diamond c=c$ but $c(x,t)$ is unbounded on some time interval
$[0,T]$. 

\section{Strong solutions}

With the same assumptions on $K$ and $\omega_{0}$ as in Section
4, we show that there is a weak solution to \eqref{eq:M_SDE1} by
using the result in \citep{Veretennikov 1981} for multi-dimensional
diffusion process with bounded drifts. Moreover, under a growth condition
on $K$, we are able to show the H\"older continuity of the vector
field $K\diamond b$.

Firstly, by using the results in Zvonkin-Krylov \citep{Zvonkin and Krylov 1974}
and Veretennikov \citep{Veretennikov 1981} we deduce the following.
\begin{thm}
\label{thm:strong1} Let $B=(B_{t})_{t\geq0}$ be a $d$-dimensional
standard Brownian motion on a probability space $(\varOmega,\mathcal{F},\mathbb{P})$.
There is a unique family of stochastic processes $X(x,t)$ which is
jointly continuous in $(x,t)$ almost surely, and satisfies the stochastic
integral equations
\[
X^{i}(x,t)=x+\int_{0}^{t}\left[\int_{\R^{d}}\left(K_{j}^{i}\star\mathscr{L}(X(y,s))\right)(X(x,s))\omega_{0}^{j}(y)dy\right]ds+B_{t}^{i}
\]
for $t\in[0,T_{L}]$, and 
\begin{equation}
b^{i}(\cdot,t)=\int_{\R^{d}}\left(K_{j}^{i}\star\mathscr{L}(X(y,s))\right)\omega_{0}^{j}(y)dy\label{eq:vect-1}
\end{equation}
where $i=1,\ldots,d$, are bounded. 
\end{thm}

\begin{proof}
According to Theorem \ref{thm:main 1}, for any $t\leq T_{L}$ (and
extended it to be zero for $t>T_{L}$), there is a unique bounded
vector field $b(x,t)$ satisfying $K\diamond b=b$. Since $b$ is
bounded and Borel measurable, then there is a unique strong solution
$X(x,t)$ to the ordinary stochastic differential equation
\[
dX(x,t)=b(X(x,t),t)dt+dB_{t},\quad X(x,0)=x
\]
(see Veretennikov \citep[Theorem 1, page 388]{Veretennikov 1981}).
Hence, by Lemma \ref{lem:2}, $X(x,t)$ is the unique strong solution
to \eqref{eq:M_SDE1}. 
\end{proof}
We are going to show that the vector field \eqref{eq:vect-1} is in
fact H\"older continuous. To this end, we need the following H\"older
continuity result of the transition probability density function,
proved originally by Nash \citep{Nash 1958} and later by Aronson
\citep{Aronson 1968-main}, Fabes and Stroock \citep{fabes1986new}. We take this from \citep[Theorem II.2.12, page 340]{stroock1988diffusion}. 
\begin{lem}
\label{lem:holder} Under the same assumption as in Lemma \ref{lem:qian-sharp-bd},
there are constants $C_{H}>0$ and $\alpha\in(0,1)$ depending only
$A$ and $d$ such that 
\begin{equation}
|p_{b}(s,y,t,x)-p_{b}(s,y,\tilde{t},\tilde{x})|\leq\frac{C_{H}}{\delta^{3}}\left(\frac{|t-\tilde{t}|^{\frac{1}{2}}\vee|x-\tilde{x}|}{\delta}\right)^{\alpha}\label{eq:new9}
\end{equation}
for all $s\geq0$, $\delta^{2}\leq t-s\leq\frac{1}{\delta^{2}},$$|x-\tilde{x}|\leq\delta$,
for any $\delta>0$. 
\end{lem}

\begin{lem}
Under the same assumptions for $K$ and $\omega_{0}$ as in the previous
section, we further assume $\gamma_{1}=\gamma_{2}\equiv\gamma$, which
belongs to $[0,d)$. Suppose $|b(x,t)|\leq C_{K}$ for all $x$ and
$t$. Then $K\diamond b(x,t)$ is H\"older continuous on any compact
subset of $\mathbb{R}^{d}\times(0,T_{K}]$, where the H\"older exponent
and H\"older norm depend only on $C_{K}$. 
\end{lem}

\begin{proof}
By using Lemma \ref{lem:holder}, if $T_{K}\geq t,\tilde{t}>\delta^{2}$
and $|x-\tilde{x}|<\delta$ (for $\delta>0$ small enough), for simplicity
set
\[
H=\left|K\diamond b(x,t)-K\diamond b(\tilde{x},\tilde{t})\right|.
\]
Then 
\begin{align*}
H & =\left|\int_{\R^{d}}\left[\int_{\mathbb{R}^{d}}K(z)\omega_{0}(y)\left(p_{b}(0,y,t,x-z)-p_{b}(0,y,\tilde{t},\tilde{x}-z)\right)dz\right]dy\right|\\
 & \leq\int_{\R^{d}}\left[\int_{\mathbb{R}^{d}}|K(z)||\omega_{0}(y)||p_{b}(0,y,t,x-z)-p_{b}(0,y,\tilde{t},\tilde{x}-z)|dz\right]dy\\
 & \leq\int_{\R^{d}}\left[\int_{\{|z|<\rho\}}|K(z)||\omega_{0}(y)||p_{b}(0,y,t,x-z)-p_{b}(0,y,\tilde{t},\tilde{x}-z)|dz\right]dy\\
 & \quad+\int_{\R^{d}}\left[\int_{\{|z|\geq\rho\}}|K(z)||\omega_{0}(y)||p_{b}(0,y,t,x-z)-p_{b}(0,y,\tilde{t},\tilde{x}-z)|dz\right]dy\\
 & \leq C_{1}\left\{ \frac{C_{H}}{\delta^{d}}\left(\frac{|t-\tilde{t}|^{\frac{1}{2}}\vee|x-\tilde{x}|}{\delta}\right)^{\alpha}\right\} \frac{C_{0}}{d-\gamma}\rho^{d-\gamma}+2\frac{C_{0}C_{1}}{\rho^{\gamma}}
\end{align*}
for any $\rho>0$. Here the last inequality follows from Lemma \ref{lem:cut-small}
and Lemma \ref{lem:cut-big}. Choose $\rho>0$ such that 
\[
C_{1}\left\{ \frac{C_{H}}{\delta^{d}}\left(\frac{|t-\tilde{t}|^{\frac{1}{2}}\vee|x-\tilde{x}|}{\delta}\right)^{\alpha}\right\} \frac{C_{0}}{d-\gamma}\rho^{d-\gamma}=2\frac{C_{0}C_{1}}{\rho^{\gamma}}
\]
that is 
\[
\frac{1}{\rho}=\left\{ \frac{C_{H}}{2(d-\gamma)\delta^{d}}\left(\frac{|t-\tilde{t}|^{\frac{1}{2}}\vee|x-\tilde{x}|}{\delta}\right)^{\alpha}\right\} ^{1/d}.
\]
Then 
\[
H\leq4C_{0}C_{1}\left\{ \frac{C_{H}}{2(d-\gamma)\delta^{d}}\left(\frac{|t-\tilde{t}|^{\frac{1}{2}}\vee|x-\tilde{x}|}{\delta}\right)^{\alpha}\right\} ^{\gamma/d},
\]
which yields the claim. 
\end{proof}
\begin{cor}
\label{thm:exist-weak-sol} Under the same conditions for $\omega_{0}$
as in the previous section. Suppose the kernel $K=(K_{j}^{i})$ satisfies
the following condition: 
\begin{equation}
\left|K(x)\right|\leq C_{0}\frac{1}{|x|^{\gamma}}\quad\textrm{ for all }x\neq0\label{eq:new7-1}
\end{equation}
where $0\leq\gamma<d$ and $C_{0}>0$. Then there is a unique strong
solution $X(x,t)$ to \eqref{eq:M_SDE1}  for any $t\leq T_{L}$,  such that
\begin{equation}
b^{i}(x,t)=\int_{\R^{d}}\left(K_{j}^{i}\star\mathscr{L}(X(y,s))\right)\omega_{0}^{j}(y)dy\label{eq:main-sde}
\end{equation}
$i=1,\cdots,d$, are bounded and H\"older continuous on any compact
subset of $\mathbb{R}^{d}\times(0,T_{L}]$.
\end{cor}

\section{From SDE to PDE}

In this section we recover the PDE from the SDE \eqref{eq:M_SDE1}.
\begin{thm}
Let $K$ and $\omega_{0}$ satisfy the assumptions in Section 4. Let
$(X(x,t),B_{t})$ (where $x\in\mathbb{R}^{d}$ and $t\geq0$) be the
unique weak solution of SDE \eqref{eq:M_SDE1} on a probability space
$(\varOmega,\mathcal{F},\mathbb{P})$ for $t\in[0,T_{L}]$. Then for
any $y\in\mathbb{R}^{d}$ and $t>0$, the distribution of $X(y,t)$
has a positive and continuous density denoted by $p(0,y,t,\cdot)$.
Let $b(x,t)$ be defined by 
\[
b^{i}(x,t)=\int_{\R^{d}}\left(K_{j}^{i}\star\mathscr{L}(X(y,s))\right)\omega_{0}^{j}(y)dy
\]
for $i=1,\ldots,d$, and
\begin{equation}
\omega^{i}(x,t)=\int_{\mathbb{R}^{d}}p(0,y,t,x)\omega_{0}^{i}(y)dy\label{eq:def-w1}
\end{equation}
for any $x$ and $t\in[0,T_{L}]$. Then the pair $(b,\omega)$ is
the solution to the following non-local partial differential equation
\begin{equation}
\frac{\partial}{\partial t}\omega^{i}+b^{j}\frac{\partial}{\partial x^{j}}\omega^{i}=\frac{1}{2}\Delta\omega^{i}-\frac{\partial b^{j}}{\partial x^{j}}\omega^{i},\quad\omega^{i}(\cdot,0)=\omega_{0}^{i}\label{eq:PDE1}
\end{equation}
and 
\begin{equation}
b^{i}(x,t)=\int_{\mathbb{R}^{d}}K_{j}^{i}(x-y)\omega^{j}(y,t)dy\label{eq:PDE2}
\end{equation}
for any $x$ any $t\in[0,T_{L}]$, where $i=1,\ldots,d$.
\end{thm}

\begin{proof}
According to our construction $b(x,t)$ is the unique bounded vector
field such that $K\diamond b=b$, and $X(x,t)$ is the unique weak
solution of the SDE 
\[
dX(x,t)=b\left(X(x,t),t\right)dt+dB_{t},\quad X(x,0)=x.
\]
Thus $p(0,y,t,x)=p_{b}(0,y,t,x)$ is the transition probability density
for the diffusion with its generator $L_{b}$, hence considering $p_{b}(0,y,t,x)$
as a function of $(t,x)$, $p_{b}$ is the fundamental solution to
the forward adjoint equation 
\[
\left(\frac{\partial}{\partial t}-L_{b}^{\star}\right)p_{b}=0
\]
where $L_{b}^{\star}=\frac{1}{2}\Delta-b\cdot\nabla-\nabla\cdot b$.
Hence, according to Lemma \ref{lem:fund-sol-vortex} $\omega(x,t)$
given by \eqref{eq:def-w1} is the solution to
\[
\left(\frac{\partial}{\partial t}-L_{b}^{\star}\right)\omega=0,\quad\omega(\cdot,0)=\omega_{0},
\]
that is
\[
\frac{\partial}{\partial t}\omega+b\cdot\nabla\omega=\frac{1}{2}\Delta\omega-\left(\nabla\cdot b\right)\omega,\quad\omega(\cdot,0)=\omega_{0}.
\]
Moreover, 
\begin{align*}
b^{i}(x,t) & =\int_{\R^{d}}\left(K_{j}^{i}\star\mathscr{L}(X(y,s))\right)\omega_{0}^{j}(y)dy\\
 & =\int_{\mathbb{R}^{d}}\left(\int K_{j}^{i}\left(x-z\right)p_{b}(0,y,t,z)dz\right)\omega_{0}^{j}(y)dy\\
 & =\int_{\mathbb{R}^{d}}K_{j}^{i}(x-z)\omega^{j}(z,t)dz
\end{align*}
which completes the proof.
\end{proof}
As an example we may apply this representation theorem to the Biot-Savart
kernel $G(x)=-\frac{x}{|x|^{3}}$ on $\mathbb{R}^{3}$ so that $K_{j}^{i}=\varepsilon^{ikj}G^{k}$.
Then there is a unique weak solution to the following SDE
\begin{equation}
dX(x,t)=\left(\int_{\R^{3}}\left.\mathbb{E}\left[G(z-X(y,t))\right]\right|_{z=X(x,t)}\wedge\omega_{0}(y)dy\right)dt+\sqrt{2\nu}dB_{t}\label{eq:sde-v3}
\end{equation}
where $\omega_{0}\in L^{1}(\mathbb{R}^{3})\cap L^{\infty}(\mathbb{R}^{3})$
is the initial vorticity. In this case we define
\[
\omega(x,t)=\int_{\mathbb{R}^{3}}p(0,y,t,x)\omega_{0}(y)dy
\]
where $p(0,y,t,\cdot)$ is the probability density function of the
law of $X(y,t)$ to \eqref{eq:sde-v3} and define
\begin{align*}
u^{i}(x,t) & =\int_{\mathbb{R}^{3}}K_{j}^{i}(x-z)\omega^{j}(z,t)dz\\
 & =-\int_{\mathbb{R}^{3}}\varepsilon^{ikj}\frac{x^{k}-z^{k}}{|x-z|^{3}}\omega^{j}(z,t)dz.
\end{align*}
That is
\[
u(x,t)=-\int_{\mathbb{R}^{d}}\frac{x-z}{|x-z|^{3}}\wedge\omega(z,t)dz.
\]
By the previous theorem $(u,\omega)$ satisfies the following PDE
\[
\frac{\partial}{\partial t}\omega+u\cdot\nabla\omega=\nu\Delta\omega-\left(\nabla\cdot u\right)\omega,\quad\omega(\cdot,0)=\omega_{0}.
\]
Moreover one can verify easily that
\[
\nabla\cdot u(x,t)=0\quad\textrm{ and }\Delta u(x,t)=-\nabla\wedge\omega(x,t)
\]
where the second equation follows from the Green formula which can
be also written as
\[
\nabla\wedge\left(\nabla\wedge u-\omega\right)=0
\]
which yields that 
\[
\omega=\nabla\wedge u-\nabla f
\]
for some scalar function $f$. If one also imposes a constrain $\nabla\cdot\omega_{0}=0$
then $\omega=\nabla\wedge u$. In this case
\[
\frac{\partial}{\partial t}\omega+u\cdot\nabla\omega=\nu\Delta\omega,\quad\omega(\cdot,0)=\omega_{0}
\]
and
\[
\omega=\nabla\wedge u.
\]
Hence $(X(x,t),B_{t})$ is the probability representation of the solution
to the above vorticity equation.

\end{document}